\documentclass{amsart}

\usepackage{amsmath}
\usepackage{amsfonts}
\usepackage{amssymb}
\usepackage{mathrsfs} 
\usepackage{enumerate}

\newtheorem*{prop}{Proposition}
\newtheorem*{thm}{Theorem}
\newtheorem*{lem}{Lemma}

\newcommand{\eps}{\varepsilon}
\newcommand{\N}{\mathbb N}

\author{J. Wengenroth}
\address{Universit\"at Trier, FB IV -- Mathematik, 54286 Trier, Germany}
\email{wengenroth@uni-trier.de}
\title{Series representation in tensor products of Banach spaces}

\begin{document}
\begin{abstract} We show by a ridiculously simple argument that, for any norm $\alpha$ on the tensor product $X\otimes Y$ of vector spaces, every element $u$ of the completion
$X\hat\otimes_\alpha Y$ can be represented as a convergent series $u=\sum\limits_{n=1}^\infty x_n\otimes y_n$ of elementary tensors. 
\end{abstract}

\maketitle
Whereas all books like, e.g., \cite{DeFl,Jar,Ryan} treating topological tensor products of Banach (or locally convex) spaces include Grothendieck's  theorem \cite[chapitre I, page 51]{Gro}that every element
of the completed {\em projective} tensor product $X \hat\otimes_\pi Y$ of Banach (or Fr\'echet) spaces has a representation as a
series of elementary tensors (even with absolute convergence)
\[
u=\sum_{n=1}^\infty x_n\otimes y_n
\] 
they all remain silent concerning the obvious corresponding question for 
the injective tensor product $X \hat\otimes_\eps Y$.
 In \cite[page 46]{Ryan} it is even claimed that {\it there is no general representation of the elements of the completed tensor product} $X \hat\otimes_\varepsilon Y$. 
Using a simplistic version of a trick of Pe{\l}czy\'nski \cite{Pel} which the author learned from \cite{Bill} we will 
show that this claim is not true.

\begin{thm}
 For any norm $\alpha$ on the tensor product $X\otimes Y$ of two vector spaces, every element $u\in X\hat\otimes_\alpha Y$ of the completion 
 has a representation $u=\sum\limits_{n=1}^\infty x_n \otimes y_n$
with $x_n \in X$ and $y_n \in Y$.
\end{thm}

At least superficially, this is similar to a Schmidt representation of compact operators between Hilbert spaces.

\begin{proof}
 Writing $u$ as the limit of a fast converging sequence of elements of $X \otimes Y$ we may form a telescoping series to obtain
$u=\sum\limits_{n=1}^\infty u_n$ with $u_n\in X\otimes Y$ such that $\sum\limits_{n=1}^\infty \alpha(u_n)<\infty$. Choosing representations of $u_n$ as sums
of $k(n)$ elementary tensors and ordering these sums one after another would give us a series $s_m$ where the 
{\it subsequence} $s_{k(1)+\cdots + k(n)}$ converges to $u$.  This, of course, does  not imply the convergence of the full
series but it would be enough if we can find representations
\[
u_n=\sum_{\ell=k(n-1)+1}^{k(n)} x_\ell\otimes y_\ell \text{ such that } \alpha\left(\sum_{\ell=k(n-1)+1}^{p} x_\ell\otimes y_\ell\right) \le 2 \alpha(u_n)
\]
for all $k(n-1)<p\le k(n)$. Indeed, for $\eps>0$, we choose $n_0\in\N$ such that $\alpha(s_{k(1)+\cdots+k(n)}-u)\le \eps$ for all $n\ge n_0$. If then $m> k(1)+\cdots+k(n_0)$
and $n\in\N$ is maximal with $k(1)+\cdots+k(n)<m$, we obtain 
\[
\alpha(s_m-u)\le \eps+ \alpha\left(\sum_{\ell=k(n)+1}^m x_\ell\otimes y_\ell\right)\le \eps+ 2\alpha(u_{n+1})
\]
which implies convergence of the series. The result thus follows from the next lemma.
\end{proof}

\begin{lem}
 For any norm $\alpha$ on the tensor product $X\otimes Y$ of two vector spaces, every element $u\in X \otimes Y$ 
 has a representation $u=\sum\limits_{\ell=1}^N x_\ell \otimes y_\ell$
with $x_\ell \in X$ and $y_\ell \in Y$ such that
\[
\alpha\left(\sum_{\ell=1}^p x_\ell\otimes y_\ell\right)\le 2 \alpha(u) \text{ for all $1\le p\le N$.}
\]
\end{lem}

For the projective norm $\alpha=\pi$ with two normed spaces, the lemma obviously holds for a  representation $u=\sum\limits_{\ell=1}^N x_\ell \otimes y_\ell$ such that
$\sum\limits_{\ell=1}^N \|x_\ell\| \|y_\ell\|\le 2\pi(u)$. This gives a very elementary proof of Grothendieck's theorem including
absolute convergence of the series.

\begin{proof} We start with an arbitrary representation $u=\sum\limits_{k=1}^m e_k\otimes f_k$ with $e_k\in X$ and $f_k\in Y$. The simplistic version of Pe{\l}czy\'nski's trick is to write, for suitably chosen $n\in\N$,
\begin{align*}
u &= \frac{1}{n} \sum_{j=1}^n u \\
& = \frac 1n e_1\otimes f_1 + \frac 1n e_2\otimes f_2 + \cdots + \frac 1n e_m\otimes f_m \\
& + \frac 1n e_1\otimes f_1 + \frac 1n e_2\otimes f_2 + \cdots + \frac 1n e_m\otimes f_m \\
&\quad \vdots \\
& + \frac 1n e_1\otimes f_1 + \frac 1n e_2\otimes f_2 + \cdots + \frac 1n e_m\otimes f_m \\
\end{align*}
where each of the $n$ rows adds to $\frac 1n u$. We  put $N=nm$ and rename the elementary tensors row after row as $x_\ell\otimes y_\ell$. This gives a representation of $u$ with the
desired property. Indeed, $p\in \{1,\ldots,N\}$ can uniquely be written as $p=qm+r$ with $q,r\in\N_0$ such that $0\le r<m$. The sum of the first 
$p$ terms $x_\ell\otimes y_\ell$ contains $q$ full rows of the
matrix above (which sum to $\frac qn u$) and the first $r$ terms of the following row. This shows
\begin{align*}
\alpha\left(\sum_{\ell=1}^p x_\ell\otimes y_\ell\right)  &= \alpha\left(\frac qn u +\sum_{k=1}^r \frac 1n e_k\otimes f_k\right) 
\le \frac qn \alpha(u) + \frac 1n \sum_{k=1}^r \alpha(e_k\otimes f_k)\\
&\le \alpha(u) +\frac mn \max\{\alpha(e_k\otimes f_k): 1\le k\le m\}.
\end{align*}
 For $\alpha(u)\neq 0$ and $n$ large enough,  this last expression is $\le 2\alpha(u)$. If $\alpha(u)=0$ we have $u=0$ and hence
can represent $u$ as an empty sum of elementary tensors. (If $\alpha$ is only a seminorm and $\alpha(u)=0$ we can achieve a representation
with $\alpha\left(\sum_{k=1}^p x_k\otimes y_k\right) \le \eps$ for any given $\eps$).
\end{proof}

We do not claim that the series representation in the theorem is of practical relevance, in particular, because there is no control about the rank of the approximating partial sums.
Nevertheless, we have the following application which is implied by another result of Grothendieck \cite[chapitre I, page 90]{Gro} about the injective tensor product, namely that $C(K)\hat\otimes_\eps Y$
is isomorphic to the Banach space $C(K,Y)$ of Banach space valued functions on a compact topological space $K$ with the uniform norm (the natural isomorphism 
maps elementary tensors $f\otimes y$ to the function $t\mapsto f(t)y$):

\begin{quote} \it
Every continuous function $F\in C(K,Y)$ on a compact topological space with values in a Banach space $Y$ has a uniformly convergent representation 
\[
F(t)=\sum\limits_{n=1}^\infty f_n(t)y_n
\]
with $f_n\in C(K)$ and $y_n\in Y$. 
\end{quote}

It is not difficult to extend the theorem to the case of metrizable locally convex topologies on $X\otimes Y$ (the representation of $u_n$ from the lemma is chosen according to the $n$-th
seminorm $\alpha_n$ on $X\otimes Y$ for an increasing sequence of seminorms $\alpha_n$ generating the topology, the final remark in the proof of the lemma concers the possibility $\alpha_n(u_n)=0$
and then should be applied, e.g., to $\eps=1/n^2$). The arguments which never used any properties of tensor products show in fact the following result.

\begin{prop}
Let $A$ be subset of a metrizable locally convex space $X$ with dense linear span. Every element $u$ of the completion $\hat X$ of $X$ has a representation
\[
u=\sum_{n=1}^\infty \lambda_n a_n \text{ with $\lambda_n\in\mathbb K$ and $a_n\in A$.}
\]
\end{prop}


The question answered in this note is somewhat naive. A perhaps better question asks for a representation
as an {\it unconditionally} convergent series. For two Hilbert spaces $X,Y$, the Schmidt representation of compact operators
$T\in \mathscr K(X',Y)=X\hat\otimes_\eps Y$ indeed yields an unconditional representation. More generally, for the completed injective tensor product
of two Banach spaces, it is sufficient that one of them has an unconditional basis but, for the general case, the present author does not know the answer. It would be enough to have in the lemma above the
stronger condition $\alpha\left(\sum\limits_{k\in J} x_k\otimes y_k\right)\le c\alpha(u)$ for a constant $c$ and all subsets $J\subseteq \{1,\ldots,N\}$.

\bibliographystyle{amsalpha}
\bibliography{Tensor-literature}

\end{document}